\documentclass{article}
\usepackage{ascmac,amsmath,amssymb,amsthm}
\newtheorem{thm}{\bf{Theorem}}

\newtheorem{lem}{\bf{Lemma}}

\allowdisplaybreaks[4]
\begin{document}
\title{An explicit upper bound of the argument of Dirichlet $L$-functions on the generalized Riemann hypothesis}
\author{Takahiro Wakasa }
\date{Graduate School of Mathematics, Nagoya University, Chikusa-ku, Nagoya 464-8602, Japan;\\ e-mail: d11003j@math.nagoya-u.ac.jp} 
\maketitle
\begin{abstract}
We prove an explicit upper bound of the function $S(t,\chi)$, defined by the argument of Dirichlet $L$-functions attached to a primitive Dirichlet character $\chi\pmod {q>1}$. An explicit upper bound of the function $S(t)$, defined by the argument of the Riemann zeta-function, have been obtained by A. Fujii \cite{Fujii2}. Our result is obtained by applying the idea of Fujii's result on $S(t)$. The constant part of the explicit upper bound of $S(t,\chi)$ in this paper does not depend on $\chi$.
\end{abstract}
\section{Introduction}
We consider the argument of Dirichlet $L$-functions. Let $L(s,\chi)$ be the Dirichlet $L$-function, where $s=\sigma+it$ is a complex variable, assosiated with a primitive Dirichlet character $\chi\pmod {q>1}$. We denote the non-trivial zeros of $L(s,\chi)$ by $\rho(\chi)=\beta(\chi)+i\gamma(\chi)$, where $\beta(\chi)$ and $\gamma(\chi)$ are real numbers. Then, when $t\neq\gamma(\chi)$, we define 
\begin{align*}
S(t,\chi)=\frac{1}{\pi}\arg L\left(\frac{1}{2}+it,\chi\right).
\end{align*}
This is given by continuous variation along the straight line $s=\sigma +it$, as $\sigma$ varies from $+\infty$ to $\frac{1}{2}$, starting with the value zero. Also, when $t=\gamma(\chi)$, we define
\begin{align*}
S(t,\chi)=\frac{1}{2}\{S(t+0,\chi)+S(t-0,\chi)\}.
\end{align*}

In Selberg \cite{Selberg2}, it is known that 
\begin{align*}
S(t,\chi)=O(\log q(t+1))
\end{align*}
and under the generalized Riemann hypothesis (GRH)
\begin{align*}
S(t,\chi)=O\left(\frac{\log q(t+1)}{\log\log q(t+3)}\right).
\end{align*}
The purpose of the present article is to prove the following result.
\begin{thm} \label{th}
Assuming GRH. Then,
\begin{align*}
|S(t,\chi)|<0.804\cdot\frac{\log q(t+1)}{\log\log q(t+3)}+O\left(\frac{\log q(t+3)}{(\log\log q(t+3))^2}\right).
\end{align*}
\end{thm}
The constant $0.804$ obviously does not depend on $\chi$. Also, the implied constant of the error term does not depend on $q$. The details of the argument concerning error terms can be seen in the proof of this theorem. However, our result does not include the case of the function $S(t)$ which is defined by the argument of the Riemann zeta-function since we assume $q>1$. An explicit upper bound of the function $S(t)$ is obtained by A. Fujii \cite{Fujii2}, where the value is 0.83. 

The basic policy of the proof of this theorem is based on A. Fujii \cite{Fujii2}. In the proof, $S(t,\chi)$ is seperated by three parts $M_1$, $M_2$ and $M_3$. Fujii's idea of \cite{Fujii2} is applied to all parts. But we need Lemma \ref{lem1}, which is an explicit formula for $\frac{L'}{L}(s,\chi)$. This lemma is an analogue of Selberg's result.

To prove our result, we introduce some notations and prove the aforementioned Lemma \ref{lem1} in Section 2.
\section{Some notations and a lemma}
Here we introduce the following notations.

Let $s=\sigma+it$. We suppose that $\sigma\geq \frac{1}{2}$ and $t\geq2$. Let $x$ be a positive number satisfying $4\leq x\leq t^2$. Also, we put
\begin{align*}
\sigma_1=\frac{1}{2}+\frac{1}{\log x}
\end{align*}
and
\begin{align*}
\Lambda_x(n)=\left\{ \begin{array}{ll}
\Lambda(n) &~~ {\rm for}~~1\leq n\leq x,  \\
\Lambda(n)\frac{\log\frac{x^2}{n}}{\log x} &~~{\rm for} ~~x\leq n\leq x^2,\\
\end{array} \right.
\end{align*}
with
\begin{align*}
\Lambda(n)=\left\{ \begin{array}{ll}
\log p &~~ {\rm if}~n=p^k~{\rm with~ a~ prime}~p~{\rm and~ an~ integer}~k\geq 1, \\
0&~~{\rm otherwise}.\\
\end{array} \right.
\end{align*}
Using these notations, we prove the following lemma.
\begin{lem} \label{lem1}
Assume the GRH. Let $t\geq2$ and $x>0$ such that $4\leq x \leq t^2$. Then for $\sigma\geq \sigma_1=\frac{1}{2}+\frac{1}{\log x}$ we have
\begin{align*}
\frac{L'}{L}(\sigma+it,\chi )&=-\sum_{n<x^2}\frac{\Lambda_x(n)}{n^{\sigma+it}}\chi(n)-\frac{x^{\frac{1}{2}-\sigma}\left(1+x^{\frac{1}{2}-\sigma}\right)\omega}{{1-\frac{1}{e}\left(1+\frac{1}{e}\right)\omega'}}\Re\left(\sum_{n<x^2}\frac{\Lambda_x(n)}{n^{\sigma_1+it}}\chi(n)\right)\\
&~~~+\frac{x^{\frac{1}{2}-\sigma}\left(1+x^{\frac{1}{2}-\sigma}\right)\omega}{{1-\frac{1}{e}\left(1+\frac{1}{e}\right)\omega'}}\cdot\frac{1}{2}\log q(t+1)+O(x^{\frac{1}{2}-\sigma}),
\end{align*}  
where $|\omega|\leq 1$ and $-1\leq \omega' \leq 1$.
\end{lem}
This is an analogue of Lemma 2 of A. Fujii \cite{Fujii2}.
\begin{lem} \label{lem2}
Let $a=0$ if $\chi(-1)=1$, and $a=1$ if $\chi(-1)=-1$. Then, for $x>1$, $s\neq -2q-a~~(q=0,1,2,\cdots)$ and $s\neq \rho(\chi)$, we have 
\begin{align*}
\frac{L'}{L}(s,\chi)&=-\sum_{n<x^2}\frac{\Lambda_x(n)}{n^{s}}\chi(n)+\frac{1}{\log x}\sum_{q=0}^{\infty}\frac{x^{-2q-a-s}-x^{-2(2q+a+s)}}{(2q+a+s)^2}\\&~~~+\frac{1}{\log x}\sum_{\rho}\frac{x^{\rho-s}-x^{2(\rho-s)}}{(s-\rho)^2}.
\end{align*}
\end{lem}
Lemma \ref{lem2} is similar to Lemma 15 of Selberg \cite{Selberg2}. We write here only a sketch of the proof of Lemma \ref{lem2}.

If $a=\max(1,\sigma)$, we have
\begin{align*}
\sum_{n<x^2}\frac{\Lambda_x(n)}{n^s}\chi(n)=\frac{1}{2\pi i\log x}\int_{a-\infty i}^{a+\infty i}\frac{x^{z-s}-x^{2(z-s)}}{(z-s)^2}\cdot\frac{L'}{L}(z,\chi)dz.
\end{align*}
We consider residues which we encounter when we move the path of integration to the left. At the point $z=s$, the residue is $-(\log x)\frac{L'}{L}(s,\chi)$. At the zeros $-2q-a~~(q=0,1,2,\cdots)$, the residues are $\frac{x^{-2q-a-s}-x^{-2(2q+a+s)}}{(2q+a+s)^2}$. At the zeros $s=\rho$ of $L(s,\chi)$, the residues are $\frac{x^{\rho-s}-x^{2(\rho-s)}}{(s-\rho)^2}$. Thus, we obtain Lemma \ref{lem2}.\\
\textit{Proof of Lemma \ref{lem1}}.
Assume the GRH. In Lemma \ref{lem2}, since for $\sigma \geq \sigma_1=\frac{1}{2}+\frac{1}{\log x}$
\begin{align*}
\left|\frac{1}{\log x}\sum_{\rho}\frac{x^{\rho-s}-x^{2(\rho-s)}}{(s-\rho)^2}\right| 
&= \frac{1}{\log x}\left|\sum_{\gamma}\frac{x^{\left(\frac{1}{2}-\sigma\right)}\left(x^{(\gamma-t)i}-x^{\left(\frac{1}{2}-\sigma\right)+2(\gamma-t)i}\right)}{\left(\sigma-\frac{1}{2}\right)^2+(t-\gamma)^2}\right|\\
&\leq \frac{x^{\frac{1}{2}-\sigma}}{\log x}\sum_{\gamma}\frac{1+x^{\frac{1}{2}-\sigma}}{\left(\sigma-\frac{1}{2}\right)^2+(t-\gamma)^2}\\
&\leq x^{\frac{1}{2}-\sigma}\left(1+x^{\frac{1}{2}-\sigma}\right)\sum_{\gamma}\frac{\sigma_1-\frac{1}{2}}{\left(\sigma_1-\frac{1}{2}\right)^2+(t-\gamma)^2},
\end{align*} 
we have 
\begin{align*}
\frac{1}{\log x}\sum_{\rho}\frac{x^{\rho-s}-x^{2(\rho-s)}}{(s-\rho)^2}= x^{\frac{1}{2}-\sigma}\left(1+x^{\frac{1}{2}-\sigma}\right)\omega\sum_{\gamma}\frac{\sigma_1-\frac{1}{2}}{\left(\sigma_1-\frac{1}{2}\right)^2+(t-\gamma)^2},
\end{align*}
where $|\omega|\leq 1$. Hence by Lemma \ref{lem2}, we have for $\sigma\geq\sigma_1$ 
\begin{align}
\frac{L'}{L}(\sigma+it ,\chi)&=-\sum_{n<x^2}\frac{\Lambda_x(n)}{n^{\sigma+it}}\chi(n)+O\left(\frac{x^{\frac{1}{2}-\sigma}}{\log x}\right)\nonumber\\
&~~~+x^{\frac{1}{2}-\sigma}\left(1+x^{\frac{1}{2}-\sigma}\right) \omega\sum_{\gamma}\frac{\sigma_1-\frac{1}{2}}{\left(\sigma_1-\frac{1}{2}\right)^2+(t-\gamma)^2}.\label{1}
\end{align}
In particular, by $x^{\frac{1}{2}-\sigma}\leq x^{-\frac{1}{\log x}}=\frac{1}{e}$ we get for $\sigma\geq\sigma_1$
\begin{align}
\Re\frac{L'}{L}(\sigma_1+it ,\chi)&=-\Re\left(\sum_{n<x^2}\frac{\Lambda_x(n)}{n^{\sigma_1+it}}\chi(n)\right)+O\left(\frac{1}{\log x}\right)\nonumber\\
&~~~+\frac{1}{e}\left(1+\frac{1}{e}\right)\omega'\sum_{\gamma}\frac{\sigma_1-\frac{1}{2}}{\left(\sigma_1-\frac{1}{2}\right)^2+(t-\gamma)^2},\label{2}
\end{align}
where $-1\leq\omega'\leq 1$.

Here, since by p. 46 of Selberg \cite{Selberg2}
\begin{align*}
\Re\frac{L'}{L}(s,\chi)=\Re\left(-\frac{1}{2}\log \frac{q}{\pi}-\frac{1}{2}\log \left(\frac{s+a}{2}\right)\right)+\sum_{\gamma}\frac{\sigma-\frac{1}{2}}{\left(\sigma-\frac{1}{2}\right)^2+(t-\gamma)^2}+O(1),
\end{align*}
we get for $t\geq2$
\begin{align}
\Re\frac{L'}{L}(\sigma_1+it,\chi)=-\frac{1}{2}\log q(t+1)+\sum_{\gamma}\frac{\sigma_1-\frac{1}{2}}{\left(\sigma_1-\frac{1}{2}\right)^2+(t-\gamma)^2}+O(1)\label{3}.
\end{align}
By (\ref{2}) and (\ref{3}) we have 
\begin{align*}
&\left(1-\frac{1}{e}\left(1+\frac{1}{e}\right)\omega'\right)\sum_{\gamma}\frac{\sigma_1-\frac{1}{2}}{\left(\sigma_1-\frac{1}{2}\right)^2+(t-\gamma)^2}\\
&~~~~~~~~~=-\Re\left(\sum_{n<x^2}\frac{\Lambda_x(n)}{n^{\sigma_1+it}}\chi(n)\right)+\frac{1}{2}\log q(t+1)+O\left(\frac{1}{\log x}\right)+O(1).
\end{align*}
Inserting the above inequality to (\ref{1}), we obtain Lemma \ref{lem1}.\\
\qed
\section{Proof of Theorem 1}
The quantity $S(t,\chi)$ is separated into the following three parts.
\begin{align*}
S(t,\chi)&=-\frac{1}{\pi}\int_{\frac{1}{2}}^{\infty}\Im\frac{L'}{L}(\sigma+it,\chi)d\sigma\\
&=-\frac{1}{\pi}\Biggl\{\Im\int_{\sigma_1}^{\infty}\frac{L'}{L}(\sigma+it,\chi)d\sigma+\Im\left\{\left(\sigma_1-\frac{1}{2}\right)\frac{L'}{L}(\sigma_1+it,\chi)\right\}\\
&~~~-\Im\int_{\frac{1}{2}}^{\sigma_1}\left\{\frac{L'}{L}(\sigma_1+it,\chi)-\frac{L'}{L}(\sigma+it,\chi)\right\}d\sigma\Biggr\}\\
&=-\frac{1}{\pi}\Im\left(M_1+M_2+M_3\right),
\end{align*}
say. 

First, we estimate $M_1$. By Lemma \ref{lem1} we have 
\begin{align}
M_1&=\int_{\sigma_1}^{\infty}\Biggl\{-\sum_{n<x^2}\frac{\Lambda_n(x)}{n^{\sigma+it}}\chi(n)-\frac{x^{\frac{1}{2}-\sigma}\left(1+x^{\frac{1}{2}-\sigma}\right)\omega}{1-\frac{1}{e}\left(1+\frac{1}{e}\right)\omega'}\Re\left(\sum_{n<x^2}\frac{\Lambda_n(x)}{n^{\sigma_1+it}}\chi(n)\right) \nonumber\\
&~~~+\frac{x^{\frac{1}{2}-\sigma}\left(1+x^{\frac{1}{2}-\sigma}\right)\omega}{1-\frac{1}{e}\left(1+\frac{1}{e}\right)\omega'}\cdot\frac{1}{2}\log q(t+1)+O\left(x^{\frac{1}{2}-\sigma}\right)\Biggr\}d\sigma \nonumber\\
&=-\int_{\sigma_1}^{\infty}\sum_{n<x^2}\frac{\Lambda_n(x)}{n^{\sigma+it}}\chi(n)d\sigma+\eta_1(t)=-\sum_{n<x^2}\frac{\Lambda_n(x)}{n^{\sigma_1+it}\log n }\chi(n)+\eta_1(t),\label{p1}
\end{align}
say. Here,
\begin{align}
|\eta_1(t)|&=\left|\int_{\sigma_1}^{\infty} \frac{\Re \omega\cdot\left(1+x^{\frac{1}{2}-\sigma}\right)x^{\frac{1}{2}-\sigma}}{1-\frac{1}{e}\left(1+\frac{1}{e}\right)\omega'}d\sigma \right|\cdot \left|\Re\left(\sum_{n<x^2}\frac{\Lambda_n(x)}{n^{\sigma_1+it}}\chi(n)\right)  -\frac{1}{2}\log q(t+1)\right|\nonumber\\
&~~~+O\left(\int_{\sigma_1}^{\infty}x^{\frac{1}{2}-\sigma}d\sigma\right)\nonumber\\
&\leq\frac{1}{1-\frac{1}{e}\left(1+\frac{1}{e}\right)}\left|\Re\left(\sum_{n<x^2}\frac{\Lambda_n(x)}{n^{\sigma_1+it}}\chi(n)\right)  -\frac{1}{2}\log q(t+1)\right|\nonumber\\
&~~~\times \int_{\sigma_1}^{\infty}x^{\frac{1}{2}-\sigma}\left(1+x^{\frac{1}{2}-\sigma}\right)d\sigma+O\left(\int_{\sigma_1}^{\infty}x^{\frac{1}{2}-\sigma}d\sigma\right)\nonumber\\
&\leq \frac{\left(\frac{1}{e}+\frac{1}{2e^2}\right)}{1-\frac{1}{e}\left(1+\frac{1}{e}\right)}\cdot\frac{1}{2}\cdot\frac{\log q(t+1)}{\log x}+O\left(\frac{1}{\log x}\left|\sum_{n<x^2}\frac{\Lambda_n(x)}{n^{\sigma_1+it}}\chi(n)\right|\right),\label{p2}
\end{align}
say.

Next, applying Lemma \ref{lem1} to $M_2$, we get
\begin{align}
|M_2|&=\Biggl|\frac{1}{\log x}\Biggr\{-\sum_{n<x^2}\frac{\Lambda_n(x)}{n^{\sigma_1+it}}\chi(n)-\frac{x^{\frac{1}{2}-\sigma_1}\left(1+x^{\frac{1}{2}-\sigma_1}\right)\omega}{1-\frac{1}{e}\left(1+\frac{1}{e}\right)\omega'}\Re\left(\sum_{n<x^2}\frac{\Lambda_n(x)}{n^{\sigma_1+it}}\chi(n)\right)\nonumber\\
&~~~+\frac{x^{\frac{1}{2}-\sigma_1}\left(1+x^{\frac{1}{2}-\sigma_1}\right)\omega}{1-\frac{1}{e}\left(1+\frac{1}{e}\right)\omega'}\cdot\frac{1}{2}\log q(t+1)+O\left(x^{\frac{1}{2}-\sigma_1}\right)\Biggr\}\Biggr|\nonumber\\
&\leq\frac{\left(\frac{1}{e}+\frac{1}{e^2}\right)}{1-\frac{1}{e}\left(1+\frac{1}{e}\right)}\cdot\frac{1}{2}\cdot\frac{\log q(t+1)}{\log x}+O\left(\frac{1}{\log x}\left|\sum_{n<x^2}\frac{\Lambda_n(x)}{n^{\sigma_1+it}}\chi(n)\right|\right),\label{p3}
\end{align} 
say. 

Next we estimate $M_3$. By Lemma 16 of Selberg \cite{Selberg2} we get
\begin{align*}
|\Im (M_3)|&=\left|\int_{\frac{1}{2}}^{\sigma_1}\Im\left\{\sum_{\rho}\frac{1}{\sigma_1+it-\rho}-\sum_{\rho}\frac{1}{\sigma+it-\rho}+O(1)\right\}d\sigma\right|\\
&\leq \left|\int_{\frac{1}{2}}^{\sigma_1} \sum_{\gamma}\frac{(t-\gamma)\left\{\left(\sigma-\frac{1}{2}\right)^2-\left(\sigma_1-\frac{1}{2}\right)^2\right\}} {\left\{\left(\sigma_1-\frac{1}{2}\right)^2+(t-\gamma)^2\right\}\left\{\left(\sigma-\frac{1}{2}\right)^2+(t-\gamma)^2\right\}}d\sigma \right| +O\left(\frac{1}{\log x}\right)\\
&=N(\gamma) +O\left(\frac{1}{\log x}\right),
\end{align*}
say. If $t=\gamma$, we see $N(\gamma) =0$ easily. If $t\neq\gamma$, we have
\begin{align*}
N(\gamma)&<\sum_{\gamma}\frac{\left(\sigma_1-\frac{1}{2}\right)^2} {\left(\sigma_1-\frac{1}{2}\right)^2+(t-\gamma)^2}\int_{\frac{1}{2}}^{\sigma_1}\frac{|t-\gamma|} {\left(\sigma-\frac{1}{2}\right)^2+(t-\gamma)^2}d\sigma\\
&\leq\sum_{\gamma}\frac{\left(\sigma_1-\frac{1}{2}\right)^2} {\left(\sigma_1-\frac{1}{2}\right)^2+(t-\gamma)^2}\int_{\frac{1}{2}}^{\infty}\frac{|t-\gamma|} {\left(\sigma-\frac{1}{2}\right)^2+(t-\gamma)^2}d\sigma\\
&\leq\frac{\pi}{2\log x}\sum_{\gamma}\frac{\sigma_1-\frac{1}{2}} {\left(\sigma_1-\frac{1}{2}\right)^2+(t-\gamma)^2}
\end{align*}  
since $\sigma<\sigma_1$ for $M_3$. 

Here, by $(\ref{2})$ and $(\ref{3})$ we get
\begin{align*}
\sum_{\gamma}\frac{\sigma_1-\frac{1}{2}} {\left(\sigma_1-\frac{1}{2}\right)^2+(t-\gamma)^2}&=\frac{1}{1-\frac{1}{e}\left(1+\frac{1}{e}\right)\omega'}\cdot\frac{1}{2}\log q(t+1)\\
&~~~+O\left(\left|\sum_{n<x^2}\frac{\Lambda_n(x)}{n^{\sigma_1+it}}\chi(n)\right|\right)+O\left(\frac{1}{(\log x)^2}\right).
\end{align*}
So, 
\begin{align*}
N(\gamma)&=\frac{\pi}{4}\cdot\frac{1}{1-\frac{1}{e}\left(1+\frac{1}{e}\right)\omega'}\cdot\frac{1}{\log x}\cdot\log q(t+1)\\
&~~~+O\left(\frac{1}{\log x}\left|\sum_{n<x^2}\frac{\Lambda_n(x)}{n^{\sigma_1+it}}\chi(n)\right|\right)+O\left(\frac{1}{(\log x)^3}\right).
\end{align*}
Hence we have
\begin{align}
|\Im(M_3)|&\leq \frac{\pi}{4}\cdot\frac{1}{1-\frac{1}{e}\left(1+\frac{1}{e}\right)\omega'}\cdot\frac{1}{\log x}\cdot\log q(t+1)\nonumber\\
&~~~+O\left(\frac{1}{\log x}\left|\sum_{n<x^2}\frac{\Lambda_n(x)}{n^{\sigma_1+it}}\chi(n)\right|\right)+O\left(\frac{1}{\log x}\right)\nonumber\\
&=\eta_4(t)+O\left(\frac{1}{\log x}\left|\sum_{n<x^2}\frac{\Lambda_n(x)}{n^{\sigma_1+it}}\chi(n)\right|\right)+O\left(\frac{1}{\log x}\right),\label{p4}
\end{align}
say.

Finally, we estimate the sums on right-hand sides of (\ref{p1}), (\ref{p2}), (\ref{p3}) and (\ref{p4}). By definition of $\Lambda_x(n)$ we have
\begin{align*}
\left|\sum_{n<x^2}\frac{\Lambda_x(n)}{n^{\sigma_1+it}}\chi(n)\right|&\leq \sum_{n<x}\frac{\Lambda(n)}{n^\frac{1}{2}}+\sum_{x\leq n\leq x^2}\frac{\Lambda(n)\log\frac{x^2}{n}}{n^\frac{1}{2}}\cdot\frac{1}{\log x} \ll \frac{x}{\log x}.
\end{align*}
Similarly,
\begin{align*}
\left|\sum_{n<x^2}\frac{\Lambda_x(n)}{n^{\sigma_1+it}\log n}\chi(n)\right|\ll \frac{x}{(\log x)^2}.
\end{align*}
So, we see
\begin{align*}
|M_1|\leq \frac{\left(\frac{1}{e}+\frac{1}{2e^2}\right)}{1-\frac{1}{e}\left(1+\frac{1}{e}\right)}\cdot\frac{1}{2}\cdot\frac{\log q(t+1)}{\log x}+O\left(\frac{x}{(\log x)^2}\right),
\end{align*}
\begin{align*}
|M_2|\leq\frac{\left(\frac{1}{e}+\frac{1}{e^2}\right)}{1-\frac{1}{e}\left(1+\frac{1}{e}\right)}\cdot\frac{1}{2}\cdot\frac{\log q(t+1)}{\log x}+O\left(\frac{x}{(\log x)^2}\right),
\end{align*}
and
\begin{align*}
|M_3|\leq \eta_4(t)+O\left(\frac{x}{(\log x)^2}\right).
\end{align*}
For $\eta_1(t)$, $\eta_2(t)$, $\eta_3(t)$ and $\eta_4(t)$, taking $x=\log q(t+3)\sqrt{\log q(t+3)}$ we have
\begin{align*}
|S(t,\chi)|&<\frac{1}{\pi}\cdot\frac{1}{1-\frac{1}{e}\left(1+\frac{1}{e}\right)} \left\{ \frac{\left(\frac{1}{e}+\frac{1}{2e^2}\right)}{2}+\frac{\left(\frac{1}{e}+\frac{1}{e^2}\right)}{2}+\frac{\pi}{4}\right\}\frac{\log q(t+1)}{\log x}\\&~~~+O\left(\frac{x}{(\log x)^2}\right)\\
&=0.803986\cdots\frac{\log q(t+1)}{\log \log q(t+3)}+O\left(\frac{\log q(t+3)}{(\log \log q(t+3))^2}\right).
\end{align*}

Therefore we obtain the theorem.
\qed\\~\\

{\bf Acknowledgments}

I thank Prof. Kohji Matsumoto for his advice and patience during the preparation of this paper. I also thank Prof. Giuseppe Molteni, Prof. Yumiko Umegaki and Dr. Ryo Tanaka, who gave many important advice.

\end{document}